\documentclass[12pt]{amsart}

\usepackage{amsmath}
\usepackage{amssymb}
\begin{document}
 \newtheorem{thm}{Theorem}[section]
 \newtheorem{cor}[thm]{Corollary}
 \newtheorem{lem}[thm]{Lemma}
 \newtheorem{prop}[thm]{Proposition}
 \theoremstyle{definition}
 \newtheorem{defn}[thm]{Definition}
 \theoremstyle{remark}
 \newtheorem{rem}[thm]{Remark}
 \newtheorem*{ex}{Example}
 \numberwithin{equation}{section}

\def\max{\operatorname{min}}
\def\card{\operatorname{card}}
\def\supp{\operatorname{supp}}
\def\dist{\operatorname{dist}}
\def\prob{\operatorname{prob}}
\def\diam{\operatorname{diam}}
\def\cl{\operatorname{cl}}
\def\log{\operatorname{log}}
\def\dim{\operatorname{dim}}
\def\Lip{\operatorname{Lip}}
\def\liminf{\operatorname{lim\,\,inf}}
\def\limsup{\operatorname{lim\,\,sup}}
\def\({\left(}
\def\){\right)}
\def\[{\left[}
\def\]{\right]}
\def\max{\operatorname{max}}
\def\Li{\operatorname{Li}}
\def\Ls{\operatorname{Ls}}
\def\Lt{\operatorname{Lt}}
\def\colon{{:}\;}

%\numberwithin{equation}{section}
%\theoremstyle{plain}
%\newtheorem{thm}{Theorem}[section]
%\endlocaldefs

\title{Feller processes on non--locally compact spaces}
\author{Tomasz Szarek}

\date\today
\maketitle

\begin{abstract}
We consider Feller processes on a complete separable metric space $X$ satisfying the ergodic condition of the form
$$
\operatornamewithlimits{lim\,\,sup}_{n\rightarrow\infty}\( \frac{1}{n} \sum_{i=1}^n P^i (x, O)\)>0\qquad\text{for some $x\in X$},
$$
where $O$ is an arbitrary open neighbourhood of some point $z\in X$ and $P$ is a transition function. It is shown that e--chains which satisfy the above condition admit an invariant probability measure. Some results on stability of such processes are presented as well.
\end{abstract}

\section{Introduction}

The theory of Feller processes is still being developed (\cite{Borovkov, Costa, Hernandez_Lasserre1, Hernandez_Lasserre2, Jarner_Tweedie1, Jarner_Tweedie2, Lasserre, Lu_Mukherjea, Meyn_Tweedie}), although these processes have been the subject of several papers over 30 years ago (see \cite{Foguel_book, Foguel1, Foguel2, Lin, Rosenblatt, Sine1, Sine2}). In most of the literature the state space is assumed compact or at least locally compact, so that existence of an invariant measure is almost immediate. In the non--locally compact case this may be proved, in turn, if a strong form of Harris reccurence on some compact set holds (see \cite{Meyn_Tweedie}). However this condition is rather very hard to verify. It is easier to obtain ergodicity on some open sets which, unfortunately, are not precompact.  Similar difficulties occur when we would like to state the Doeblin condition (see \cite{Meyn_Tweedie}). 

It seems that the non--locally compact case has not yet been completely analyzed. 
In this note we give a contribution to this. 
The work was motivated by the need to investigate the limit behaviour of discrete Markov chains generated by iterated function systems (\cite{Barnsley, Diaconis_Freedman, Jarner_Tweedie1, Lasota_Yorke, Steinsaltz}) and stochastic differential equations on Hilbert spaces (see \cite{Da Prato_Zabczyk}). The utility of our method in proving the existence of an invariant measure for stochastic partial differential equations with an impulsive noise will be shown in \cite{Lasota_Szarek}.

\vskip5mm

Let $(X, \rho)$ be a complete and separable metric space and
let $\mathbf\Phi=(\Phi_n)_{n\ge 1}$ be a discrete--time Markov chain on $X$. By $\mathcal B(X)$ we denote the space of all Borel sets. Let $P(x, A)$ be a transition function defined for $x\in X$ and $A\in\mathcal B(X)$. {\bf Feller's property} means that the function $x\rightarrow P(x, U)$ is lower semicontinuous for all open sets $U$. Alternatively we can say that
$$
C(X)\ni f(\cdot )\rightarrow Pf (\cdot)=\int_X f(y) P(\cdot, dy)\in C(X),
$$
where $C(X)$ denotes the space of all bounded continuous functions on $X$. 

We are interested in the existence of an invariant probability measure for $\mathbf\Phi$. A measure $\mu$ is called {\bf invariant} if
$$
\mu (A)=\mu P(A)=\int_X P(x, A)\mu (dx)
$$
for $A\in\mathcal B(X)$.
\vspace{0.5cm}

Let $\mu$ be an arbitrary Borel measure. We define the {\bf support of the measure} $\mu$ by setting
$$
\supp\mu=\{x\in X: \mu(B(x, \varepsilon))>0\,\,\,\text{for every}\,\,\,\varepsilon>0\}.
$$

In order to establish the existence of an invariant measure and stability we introduce the following condition:
\vspace{0.5cm}

${\bf (\mathcal E )}$ There exists $z\in X$ such that for every open set $O$ containing $z$
$$
\operatornamewithlimits{lim\,\,sup}_{n\rightarrow\infty}\( \frac{1}{n} \sum_{i=1}^n P^i (x, O)\)>0\qquad\text{for some $x\in X$}.
\leqno(1.1)
$$

\section{Existence of invariant measures}
\begin{prop} Let $P\colon X\times \mathcal B (X)\rightarrow [0, 1]$ be a transition function for a discrete--time Markov chain $\mathbf\Phi$ and assume that condition $(\mathcal E)$ holds for some $z\in X$. If $\{P^n f\colon n\in\mathbb N\}$ is equicontinuous in $z$ for every Lipschitz continuous function $f$, then $\mathbf\Phi$ admits an invariant probability measure.
\end{prop}

\begin{proof} To finish the proof it suffices to show that for every $\varepsilon>0$ there exists a compact set $K\subset X$ such that
$$
\operatornamewithlimits{lim\,\,inf}_{n\rightarrow\infty}P^n (z, K^{\varepsilon})\ge 1-\varepsilon,
\leqno(2.1)
$$
where $K^{\varepsilon}=\{x\in X: \inf_{y\in K} \rho (x, y)<\varepsilon\}.$ This, in conjunction with Theorem 2.2 in \cite{Ethier_Kurtz}, tells us that the measures $\{P^n (z, \cdot): n\in\mathbb N\}$ are tight. Therefore the Cesaro averages are weakly precompact by the Prokhorov theorem (see \cite{Ethier_Kurtz}). Note that any weak limit of the Cesaro averages is invariant.

Assume, contrary to our claim, that (2.1) does not hold for some $\varepsilon>0$. By Ulam's lemma (see \cite{Billingsley}) there exist a sequence of compact sets $(K_i)_{i\ge 1}$ and a sequence of integers $(q_i)_{i\ge 1}$ satisfying
$$
P^{q_i} (z, K_i)>\varepsilon
$$
and
$$
\min\{\rho (x, y): x\in K_i, y\in K_j\}\ge \varepsilon/3\quad\text{for $i, j\in\mathbb N,\,\, i\neq j.$}
\leqno(2.2)
$$

We first show that for every open set $O$ containing $z$ and $j\in\mathbb N$ there exist $y\in O$ and $i\ge j$ such that
$$
P^{q_i}(y, K_i^{\varepsilon/12})<\varepsilon/2.
$$
On the contrary, suppose that there exist an open set $O'$ containing $z$ and $i_0\in\mathbb N$ such that 
$$
\inf\{ P^{q_i} (y, K_i^{\varepsilon/12}): y\in O',\,\,\, i\ge i_0\}\ge\varepsilon/2.
\leqno(2.3)
$$
Let $x\in X$ be such that condition (1.1) holds with $O'$ in place of $O$.
Let $\alpha>0$ be such that
$$
\operatornamewithlimits{lim\,\,sup}_{n\rightarrow\infty}\( \frac{1}{n} \sum_{i=1}^n P^i (x, O')\)>\alpha.
$$
 By (2.2), (2.3) and the Chapman--Kolmogorov equation we obtain
$$
\operatornamewithlimits{lim\,\,sup}_{n\rightarrow\infty}\frac{1}{n} \sum_{i=1}^n P^i \(x,  \bigcup_{j=i_0}^{N} K_j^{\varepsilon/12}\)> (N-i_0)\alpha\varepsilon/2
$$
for every $N\ge i_0$, which is impossible.

We will now define by induction a sequence of Lipschitz continuous functions $(\tilde {f}_n)_{n\ge 1}$
, a sequence of points $(y_n)_{n\ge 1}$, $y_n\to z$ as $n\to\infty$, and three 
increasing sequences of integers $(i_n)_{n\ge 1}$, $(k_n)_{n\ge 1}$, $(m_n)_{n\ge 1}$, $i_{n+1}> k_n>i_n$ for $n\in\mathbb N$, such that 
$$
\tilde {f_n}_{| K_{i_n}}=1\quad\text{and}\quad 0\le\tilde {f_n}\le \mathbf 1_{K_{i_n}^{\varepsilon/12}},
\leqno(2.4)
$$
$$
\left | P^{m_n} \(\sum_{i=1}^n\tilde {f_i}\) (z) - P^{m_n} \(\sum_{i=1}^n\tilde {f_i}\) (y_n)\right| > \varepsilon/4
\leqno(2.5)
$$
and 
$$
P^{m_n}\(u, \bigcup_{i=k_n}^{\infty} K_i^{\varepsilon/12}\)<\varepsilon/16\quad\text{for $u=z,\,\, y_n$,\,\,\,\, $n\in\mathbb N$}.
\leqno(2.6)
$$
Let $n=1$. From what has already been proved, it follows that there exist $y_1\in B(z, 1)$ and  $i_1\in\mathbb N$ such that
$$
P^{q_{i_1}}(y_1, K_{i_1}^{\varepsilon/12})<\varepsilon/2.
$$
Set $m_1=q_{i_1}$ and let $k_1 >i_1$ be such that
$$
P^{m_1}\(u, \bigcup_{i=k_1}^{\infty} K_i^{\varepsilon/12}\)<\varepsilon/16\quad\text{for $u=z,\,\, y_1$.}
$$
Let $\tilde {f_1}$ be an arbitrary Lipschitz continuous function satisfying
$$
{\tilde{f_1}}_{| K_{i_1}}=1\quad\text{and}\quad 0\le\tilde{f_1}\le \mathbf 1_{K_{i_1}^{\varepsilon/12}}.
\leqno(2.7)
$$
Thus
$$
|P^{m_1} \tilde{f_1} (z) - P^{m_1} \tilde{f_1} (y_1)|\ge P^{m_1}(z, K_{i_1})-P^{m_1}(y_1,K_{i_1}^{\varepsilon/12}) >\varepsilon/2.
$$
If $n\ge 2$ is fixed and $\tilde{f}_1, \ldots, \tilde{f}_{n-1}$; $y_1, \ldots, y_{n-1}$; $i_1, \ldots, i_{n-1}$; $k_1,\ldots, k_{n-1}$; $m_1,\ldots, m_{n-1}$ are given we choose $\sigma <n^{-1}$ such that
$$
\left|P^{m} \(\sum_{i=1}^{n-1}\tilde{f}_i\) (z)-P^{m} \(\sum_{i=1}^{n-1} \tilde{f}_i\) (y)\right|<\varepsilon/8
\leqno(2.8)
$$
for $y\in B(z, \sigma)$ and $m\in\mathbb N$. 
Similarly as in the first part, we may choose $y_{n}\in B(z, \sigma)$ and $i_{n}>k_{n-1}$ such that
$$
P^{q_{i_{n}}} (y_{n},  K_{i_{n}}^{\varepsilon/12})<\varepsilon/2.
$$
Set $m_{n}=q_{i_{n}}$ and let $\tilde  f_{n}$ be an arbitrary Lipschitz continuous function satisfying condition (2.4).
Let $k_{n} > i_{n}$ be such that
$$
P^{m_{n}}\(u, \bigcup_{i=k_{n}}^{\infty} K_i^{\varepsilon/12}\)<\varepsilon/16\quad\text{for $u=z,\,\, y_{n}$.}
$$
From this, (2.8) and the definition of $\tilde {f}_{n}$ we have
\begin{align*}
&\left|P^{m_{n}} \(\sum_{i=1}^{n}\tilde{f}_{i}\) (z)-P^{m_{n}} \( \sum_{i=1}^{n}\tilde{f}_{i}\) (y_{n})\right|\\
&\ge \left|P^{m_{n}} \tilde f_{n} (z)-P^{m_{n}} \tilde f_{n} (y_{n})\right|\\
&-\left|P^{m_{n}}\(\sum_{i=1}^{n-1} \tilde{f}_{i}\) (z)-P^{m_{n}} \(\sum_{i=1}^{n-1} \tilde{f}_{n}\) (y_{n})\right|\\
&>\varepsilon/2 - \varepsilon/8 >\varepsilon/4.
\end{align*}

We now define $f=\sum_{i=1}^{\infty} \tilde f_{i}$. By (2.2) and (2.4) $f$ is a Lipschitz continuous function and $\|f\|_{\infty}\le 1$. Finally, by (2.5) and (2.6) we have
$$
|P^{m_n} f(z) - P^{m_n} f(y_n)|>\varepsilon/8\qquad\text{for $n\in\mathbb N$}
$$
and since $y_n\rightarrow z$ as $n\rightarrow\infty$, this contradicts the assumption that $\{P^n f\colon n\in\mathbb N\}$ is equicontinuous in $z$.\end{proof}

The Markov transition function $P$ is called {\bf equicontinuous} if for $f\in C_b (X)$ the sequence of functions $\{P^n f \colon n\in  N\}$ is equicontinuous on compact sets. Recall that by $C_b (X)$ we denote the space of all bounded continuous functions with a bounded support.

A Markov chain which possesses an equicontinuous Markov transition function will be called an {\bf e--chain}.

\vskip5mm

{\bf Remark:} The concept of e--chains appears in  \cite{Jamison1, Jamison2, Rosenblatt, Sine1, Sine2}. It is, of course, clear that the condition appearing in the definition of an e--chain is equivalent to equicontinuity of $\{P^n f: n\in\mathbb N\}$, $f\in C_b (X)$, in every point $x\in X$. 

\vskip5mm
In Proposition 2.1 we assumed that equicontinuity holds for all Lipschitz continuous functions. Now we introduce a condition which allows to restrict to the case of all Lipschitz continuous functions with a bounded support.

\vskip5mm
A continuous function $V\colon X\to [0, \infty)$ is called a {\bf Lyapunov function} if
$$
\operatornamewithlimits{lim}_{\rho (x, x_0)\to\infty} V(x)=\infty
$$
for some $x_0\in X$.

\begin{thm} Let $\mathbf\Phi$ be an e--chain such that condition $(\mathcal E)$ holds and let $P\colon X\times\mathcal B (X)\rightarrow [0, 1]$ be its transition function. If there exist a Lyapunov function $V\colon X\rightarrow [0, \infty)$ and $\lambda <1$, $b<\infty$, $R<\infty$, $x_0\in X$ such that
$$
PV(x)\le \lambda V(x) + b\mathbf 1_{B(x_0, R)}(x)\qquad\text{for $x\in X$},
\leqno(2.9)
$$
then $\mathbf\Phi$ admits at least one invariant probability measure.
\end{thm}

\begin{proof} Observe that (2.9) implies that $\mathbf\Phi$ is {\bf bounded in probability}, i. e. for $x\in X$ and $\varepsilon>0$ there exists a bounded Borel set $C\subset X$ such that $P^n (x, C)\ge 1-\varepsilon$ for $n\in\mathbb N$ (see \cite{Jarner_Tweedie2}).
If we assume, contrary to our claim, that $\mathbf\Phi$ does not admit an invariant probability measure, the same conclusion as in the proof of Proposition 2.1 can be drawn for some Lipschitz continuous function with a bounded support.
\end{proof}

As an illustration of the power of Proposition 2.1, we have the following example:
\vskip5mm

{\bf Example:} (Jump process)
We consider a jump process connected with an iterated function system. Similar process on $\mathbb R^n$ was considered in \cite{PichorRudnicki}. Let $(\Omega, \mathcal F, \mathrm {Prob})$ be a probability space and let $(\tau_n)_{n\ge 0}$ be a sequence of random variables $\tau_n:\Omega\to \mathbb R_+$ with $\tau_0=0$ and such that $\Delta\tau_n=\tau_n-\tau_{n-1}$, $n\ge 1$, are independent and have the same density $\gamma e^{-\gamma t}$. Let $(S(t))_{t\ge 0}$ be a continuous semigroup on $X$. We have also given a sequence of continuous transformations $w_i:X\to X$, $i=1, \ldots, N$, and a probabilistic vector $(p_1(x), \ldots, p_N(x))$, $p_i (x)\ge 0$, $\sum_{i=1}^N p_i(x)=1$ for $x\in X$. The pair $(w_1, \ldots, w_N; p_1, \ldots, p_N)$ is called an {\bf iterated function system}.

Now we define the $X$--valued Markov chain $\Phi=(\Phi_n)_{n\ge 1}$ in the following way. We choose $x\in X$ and let $\xi_1=S(\tau_1)(x)$. We randomly select from the set $\{1, \ldots, N\}$ an integer $i_1$ and the probability that $i_1=k$ is equal to $p_k (\xi_1)$. Set $\Phi_1=w_{i_1} (\xi_1)$.

Let $\Phi_1, \ldots, \Phi_{n-1}$, $n\ge 2$, be given. Assuming that $\Delta\tau_n=\tau_{n}-\tau_{n-1}$ is independent upon $\Phi_1, \ldots, \Phi_{n-1}$, we define $\xi_n=S(\Delta\tau_n)(\Phi_{n-1})$. Further, we randomly choose $i_n$ from the set $\{1, \ldots, N\}$ in such a way that the probability of the event $\{i_n=k\}$ is equal to $p_k (\xi_n)$. Finally, we define $\Phi_n=w_{i_n} (\xi_n)$.

We will assume that there exists $r\in (0, 1)$ such that
$$
\sum_{i=1}^N p_i(x)\rho (w_i(x), w_i (y))\le r\rho (x, y)\quad\text{for $x, y\in X$.}
\leqno(2.10)
$$
Moreover, there exist $a>0$ such that
$$
\sum_{i=1}^N|p_i(x)-p_i(y)|\le a\rho (x, y)\quad\text{for $x, y\in X$}
\leqno(2.11)
$$
and $\kappa\ge 0$ such that
$$
\rho(S(t)(x), S(t)(y))\le e^{\kappa t}\rho(x, y)\quad\text{for $x, y\in X$ and $t\ge 0$}.
\leqno(2.12)
$$

We will assume that a semigroup $(S(t))_{t\ge 0}$ admits a global attractor. Recall that a compact set $\mathcal K\subset X$ is called a {\bf global attractor} if it is invariant and attracting for  $(S(t))_{t\ge 0}$, i.e. $S(t)\mathcal K=\mathcal K$ for every $t\ge 0$ and for every bounded ball $B$ and open set $U$, $\mathcal K\subset U$, there exists $t_*>0$ such that $S(t)B\subset U$ for $t\ge t_*$.

\begin{prop} 
Assume that conditions (2.10)--(2.12) hold and 
$$
r+\kappa/\gamma<1.
\leqno(2.13)
$$
If $(S(t))_{t\ge 0}$ has a global attractor, then $\Phi$ admits an invariant probability measure.
\end{prop}

\begin{proof} It is easily seen that $\Phi$ is a Markov chain. Analysis similar to that in \cite{Horbacz} (see also \cite{Lasota_Yorke, Horbacz_k}) shows that its transition function must be of the form
$$
P(x, A)=\sum_{i=1}^N\int_0^{\infty} \gamma e^{-\gamma t}p_i(S(t)(x)){\mathbf 1}_A(w_i (S(t)(x)))dt
\leqno(2.14)
$$
for $x\in X$ and $A\in\mathcal B(X).$
Then
$$
Pf(x)=\sum_{i=1}^N\int_0^{\infty} \gamma e^{-\gamma t}p_i(S(t)(x))f(w_i (S(t)(x)))dt
$$
for every $f\in C(X)$ and $x\in X$.

Let $L\ge a\gamma (\kappa-\gamma (1+r))^{-1}$ and let $f$ be a Lipschitz continuous function with the Lipschitz constant $L$. If $\|f\|_{\infty}\le 1$, then $\|Pf\|_{\infty}\le 1$ and
\begin{align*}
&|P f(x)-P f(y)|\\
&\le\sum_{i=1}^N\int_0^{\infty} \gamma e^{-\gamma t}p_i(S(t)(x))|f(w_i(S(t)(x)))-f(w_i(S(t)(y)))|dt\\
&+\sum_{i=1}^N\int_0^{\infty} \gamma e^{-\gamma t}|p_i(S(t)(x))-p_i(S(t)(y))|dt\\
&\le Lr\(\int_0^{\infty} \gamma e^{-\gamma t+\kappa t}dt\)\rho (x, y)
+a\(\int_0^{\infty} \gamma e^{-\gamma t+\kappa t}dt\)\rho (x, y)\\
&\le L\rho (x, y)\qquad\qquad\qquad\,\,\,\,\,\,\,\,\,\text{for $x, y\in X$.}
\end{align*}
From this and the fact that $P$ is linear it follows that $\{P^n f:n\in\mathbb N\}$ is equicontinuous in any $x\in X$ for an arbitrary Lipschitz continuous function $f$. Let $x_0\in X$ and set $V(x)=\rho (x, x_0)$ for $x\in X$. An easy computation shows that
$$
PV(x)\le r\gamma (\gamma-\kappa)^{-1} V(x)+N \tilde b\qquad\text{for $x\in X$},
$$
where $\tilde b=\sup_{t\ge 0, 1\le i\le N} \rho (w_i(S(t)(x_0)), x_0)<\infty$, by the fact that $(S(t))_{t\ge 0}$ has a global attractor. Set $\lambda_0=r\gamma (\gamma-\kappa)^{-1}$.
By (2.13) we have $\lambda_0<1$. Let $\lambda\in (\lambda_0, 1)$. Since $V$ is a Lyapunov function, there exists $R>0$ such that condition (2.9) holds with $b=N\tilde b$.
Hence $\Phi$ is bounded in probability (see \cite{Meyn_Tweedie}). Fix $x\in X$ and let $C\subset X$ be a bounded Borel set such that $P^n(x, C)> 1/2$.
Let $\mathcal K\subset X$ be an attractor for $(S(t))_{t\ge 0}$ and let $K=\bigcup_{i=1}^N w_i (\mathcal K)$. Since $w_i$, $i=1, \ldots, N$, are continuous, the set $K\subset X$ is compact. Further, from (2.14) and the fact that $\mathcal K$ was a global attractor, it follows that 
for every open set $U$, $K\subset U$, there exists a positive constant $\beta$ such that
$$
P(y, U)\ge \beta\quad\text{for $y\in C$.}
$$
This and the Chapman--Kolmogorov equation give
$$
\operatornamewithlimits{lim\,\,inf}_{n\rightarrow\infty}\( \frac{1}{n} \sum_{i=1}^n P^i (x, O)\)>\beta/2.
$$
Since $K$ is compact, we see that there exists $z\in K$ such that condition $(1.1)$ holds for every open neighbourhood $U$ of $z$. Thus $\Phi$ has an invariant measure by Proposition 2.1.
\end{proof}

\section{Stability results}

\begin{thm}
Let $\mathbf\Phi$ be an e--chain. Let $P\colon X\times\mathcal B (X)\rightarrow [0, 1]$ be its transition function and assume that there exists $z\in X$ such that for every open set $O$ containing $z$
$$
\operatornamewithlimits{lim\,\,inf}_{n\rightarrow\infty}P^n (x, O)>0\qquad\text{for $x\in X$}.
\leqno(3.1)
$$
Let
$$
\mathcal Z=\overline {\bigcup_{n=1}^{\infty} \supp P^n (z, \cdot)}.
$$
If there exist a Lyapunov function $V\colon X\rightarrow [0, \infty)$ and $\lambda <1$, $b<\infty$, $R<\infty$, $x_0\in X$ such that (2.9) holds,
then $\mathbf\Phi$ admits a unique invariant probability measure $\mu_*$ supported on $\mathcal Z$. Moreover
$$
\mu P^n\stackrel{\rm w}{\rightarrow} \mu_*\qquad\text{as $n\rightarrow\infty$}
$$
for every probability measure $\mu$ such that $\supp\mu\subset\mathcal Z$.
\end{thm}

\begin{proof} Since (3.1) implies (1.1), from Theorem 2.2 it follows that $\mathbf\Phi$ has an invariant probability measure, say $\mu_*$. It may be obtained (see \cite{Ethier_Kurtz, Szarek}) as any weak limit of the Cesaro averages of $(P^n (z, \cdot))_{n\ge 1}$. Therefore we may assume that $\supp\mu_*\subset\mathcal Z$.

Let us denote by $\Delta  (x_1, x_2; f; \varepsilon)$ for $x_1, x_2\in X$, $f\in C_b(X)$, $\varepsilon>0$ the set of all $\alpha\in (0, 1]$ such that there exist probability measures $\mu_1, \mu_2$ and an integer $m$ satisfying
$$
P^m (x_i, \cdot)\ge\alpha\mu_i(\cdot)\qquad\text{for $i=1, 2$},
\leqno(3.2)
$$
and
$$
\left |\int_X f(y) \mu_1 P^n (dy)-\int_X f(y) \mu_2 P^n (dy)\right |\le\varepsilon\qquad\text{for $n\in\mathbb N$.}
\leqno(3.3)
$$

We claim that $\sup \Delta (x_1, x_2; f; \varepsilon)=1$ for $x_1, x_2\in\mathcal Z$, $f\in C_b(X)$ and $\varepsilon>0$.
Fix $x_1, x_2\in \mathcal Z$, $f\in C_b(X)$ and $\varepsilon >0$. By the Chapman--Kolmogorov equation we easily obtain that
$$
\operatornamewithlimits{lim\,\,inf}_{n\rightarrow\infty}P^n (x, O_i)>0\qquad\text{for $x\in X$},
$$
where $O_i$ is an arbitrary open set containing $x_i$, $i=1, 2$. Now from the proof of Proposition 2.1 it follows that the families $\{P^n (x_i, \cdot)\colon n\in\mathbb N\}$, $i=1, 2,$ are weakly precompact (see also Theorem 2.2 in \cite{Ethier_Kurtz}).
Let $\sigma>0$ be such that
$$
|P^n f(z)-P^n f(y)|\le\varepsilon\quad\text{for $y\in B(z, \sigma)$ and $n\in\mathbb N$.}
\leqno(3.4)
$$
By (3.1) there exist $m\in\mathbb N$ and $\tilde\alpha>0$ such that
$$
P^m (x_i, B(z, \sigma))\ge\tilde\alpha\qquad\text{for $i=1, 2$.}
$$
Define
$$
\tilde\mu_i (\cdot)=\frac{P^m (x_i, B(z, \sigma)\cap \cdot)}{P^m (x_i, B(z, \sigma))}\qquad\text{for $i=1, 2$,}
\leqno(3.5)
$$
and observe that condition (3.2) is satisfied with $\tilde\mu_i$ in place of $\mu_i$ and $\tilde\alpha$ in place of $\alpha$. Moreover, from (3.4) it follows that (3.3) holds with $\tilde\mu_i$ in place of $\mu_i$. Hence $\Delta (x_1, x_2; f; \varepsilon)\neq\emptyset.$ Set $\alpha_0=\sup\Delta (x_1, x_2; f; \varepsilon)$. Suppose, contrary to our claim, that $\alpha_0<1$. Let $(\alpha_n)_{n\ge 1}$ be such that $\alpha_n\rightarrow\alpha_0$ as $n\rightarrow\infty$ and $\alpha_n\in\Delta (x_1, x_2; f; \varepsilon)$ for $n\in\mathbb N$. Let $\mu_i^n$, $i=1, 2$, and $m_n$ satisfy (3.2) with $\alpha_n$ in place of $\alpha$. Since $\{P^n (x_i, \cdot)\colon n\in\mathbb N\}$, $i=1, 2$, are tight, $\{P^{m_n} (x_i, \cdot)-\alpha_n\mu_i^n\colon n\in\mathbb N\}$, $i=1, 2,$ are weakly precompact. Therefore, without loss of generality, we may assume that $(P^{m_n} (x_i, \cdot)-\alpha_n\mu_i^n)_{n\ge 1}$, $i=1, 2$, converge to some measures $\tilde\mu_1$, $\tilde\mu_2$, respectively. Choose $y_1\in\supp\tilde\mu_1$ and $y_2\in\supp\tilde\mu_2$. From (3.1) it follows that there exist $m\in\mathbb N$ and $\gamma>0$ such that
$$
P^m (y_i, B(z, \sigma))\ge\gamma\qquad\text{for $i=1, 2$.}
$$
By Feller's property, there is $r>0$ such that
$$
P^m (y, B(z, \sigma))\ge\gamma/2\quad\text{for $y\in B(y_i, r)$, $i=1, 2$.}
$$
Set
$$
s_0=\min\{\tilde\mu_1 (B(y_1, r)), \tilde\mu_2 (B(y_2, r))\}
$$
and observe that $s_0>0$.
By the Alexandrov theorem (see \cite{Billingsley}) we may choose $k\in\mathbb N$ such that
$$
P^{m_k} (x_i, B(y_i, r))-\alpha_k \mu_i^k (B(y_i, r))>s_0/2\quad\text{for $i=1, 2$.}
$$
Let $k\in\mathbb N$ be such that
$$
\alpha_k+s_0\gamma/4>\alpha_0.
\leqno(3.6)
$$
Then by the Chapman--Kolmogorov equation (see also (3.5)) we obtain that there exist
probability measures $\hat\mu_i$ with $\supp\hat\mu_i\subset B(z, \sigma)$, $i=1, 2$, such that
$$
P^{m_k+m} (x_i, \cdot)-\alpha_k \mu_i^k P^m\ge s_0\gamma\hat\mu_i/4.
$$
Set 
$$
\mu_i=(\alpha_k + s_0\gamma/4)^{-1} (\alpha_k\mu_i^k P^{m} + s_0\gamma\hat\mu_i/4)\quad\text{for $i=1, 2$.}
$$
Since $\supp\hat\mu_i\subset B(z, \sigma)$ for $i=1, 2$, from (3.4) it follows that
$\mu_i$, $i=1, 2,$
satisfy (3.3). Finally, observe that $\mu_i$, $i=1, 2$, satisfy condition (3.2) with $m_k +m$ in place of $m$ and $\alpha=\alpha_k +s_0\gamma/4$. Hence 
$\alpha_k +s_0\gamma/4\in \Delta (x_1, x_2; f; \varepsilon)$, which contradicts the definition of $\alpha_0$, by (3.6).

We have proved that
$$
\operatornamewithlimits{lim}_{n\rightarrow\infty}\left|\int_X f(y)\mu_1P^n (dy)-\int_X f(y)\mu_2P^n (dy)\right|=0
$$
for all point measures $\mu_1, \mu_2$ supported on $\mathcal Z$ and for every $f\in C_b(X)$.
Since linear combinations of point measures are dense in the space of all measures equipped with the weak topology, the above convergence holds
for all probability measures $\mu_1, \mu_2$ supported on $\mathcal Z$ and for every $f\in C_b(X)$. 
Since $\mathbf \Phi$ is bounded in probability, the above convergence is also satisfied for every $f\in C(X)$.
From this it follows that $\mu_*$ is a unique invariant measure supported on  $\mathcal Z$ and 
$$
\mu P^n\stackrel{\rm w}{\rightarrow} \mu_*\qquad\text{as $n\rightarrow\infty$}
$$
for every probability measure $\mu$ such that $\supp\mu\subset\mathcal Z$, which finishes the proof.
\end{proof}

A point $x\in X$ is called {\bf reachable} if for every open set $O$ containing $x$
$$
\sum_{n=1}^{\infty} P^n (y, O)>0\quad\text{for every $y\in X$}.
$$
The chain $\mathbf \Phi$ is called {\bf open set irreducible} if every point is reachable.
\vskip3mm
As a consequence of Theorem 3.1 and the above definition we obtain the following theorem:

\begin{thm}
Let $\mathbf\Phi$ be an open set irreducible e--chain. Let $P\colon X\times\mathcal B (X)\rightarrow [0, 1]$ be its transition function and assume that there exists $z\in X$ such that for every open set $O$ containing $z$ condition (3.1) holds.
If there exist a Lyapunov function $V\colon X\rightarrow [0, \infty)$ and $\lambda <1$, $b<\infty$, $R<\infty$, $x_0\in X$ such that (2.9) holds,
then $\mathbf\Phi$ admits a unique invariant probability measure $\mu_*$. Moreover
$$
\mu P^n\stackrel{\rm w}{\rightarrow} \mu_*\qquad\text{as $n\rightarrow\infty$}
$$
for every probability measure $\mu$.
\end{thm}
\begin{proof} It suffices to note that
$$
\overline{\bigcup_{n=1}^{\infty} \supp P^n (z, \cdot)}=X.
$$
\end{proof}

\begin{thm}
Let $\mathbf\Phi$ be an e--chain. Let $P\colon X\times\mathcal B (X)\rightarrow [0, 1]$ be its transition function and assume that there exists $z\in X$ such that for every open set $O$ containing $z$ there exists $\alpha>0$ satisfying
$$
\operatornamewithlimits{lim\,\,inf}_{n\rightarrow\infty}P^n (x, O)\ge\alpha\quad\text{for $x\in X$}.
\leqno(3.7)
$$
If there exist a Lyapunov function $V\colon X\rightarrow [0, \infty)$ and $\lambda <1$, $b<\infty$, $R<\infty$, $x_0\in X$ such that (2.9) holds,
then $\mathbf\Phi$ admits a unique invariant probability measure $\mu_*$. Moreover
$$
\mu P^n\stackrel{\rm w}{\rightarrow} \mu_*\qquad\text{as $n\rightarrow\infty$}
\leqno(3.8)
$$
for every probability measure $\mu$.
\end{thm}
\begin{proof} The existence of an invariant measure $\mu_*$ follows from Theorem 2.2. Fix $\varepsilon>0$, $x_1, x_2\in X$ and $f\in C_b (X)$. By equicontinuity of $\{P^n f\colon n\in\mathbb N\}$ in $z\in X$, we choose $r>0$ such that
$$
|P^n f(z)-P^n f(x)|<\varepsilon/4\quad\text{for $x\in B(z, r)$ and $n\in\mathbb N$.}
\leqno(3.9)
$$
Let $\alpha>0$ be such that (3.7) holds with $O=B(z, r)$. Then by Fatou's lemma we have
$$
\operatornamewithlimits{lim\,\,inf}_{n\rightarrow\infty}\mu P^n (O)\ge\alpha
\leqno(3.10)
$$
for every probability measure $\mu$. Let $k\in\mathbb N$ be such that $4(1-\alpha/2)^k \|f\|_{\infty}\le\varepsilon$. Further, from the Lasota--Yorke theorem (see Theorem 4.1 in \cite{Lasota_Yorke}) and (3.10) it follows that there exist integers $n_1,\ldots, n_k$ and proba\-bility mea\-sures $\nu_1^i, \ldots, \nu_k^i, \mu_k^i$ such that $\supp\nu_j^i\subset O$, $j=1,\ldots, k,$ and
\begin{align*}
&P^{n_1+\ldots +n_k} (x_i, \cdot)=\frac{\alpha}{2}\nu_1^i P^{n_2+\ldots+n_k}+\frac{\alpha}{2}\(1-\frac{\alpha}{2}\)\nu_2^iP^{n_3+\ldots+n_k}\\
&+\ldots+\frac{\alpha}{2}\(1-\frac{\alpha}{2}\)^{k-1}\nu_k^i+\(1-\frac{\alpha}{2}\)^{k}\mu_k^i\quad\text{for $i=1, 2$}.
\end{align*}
Then by the Markov property we obtain 
\begin{align*}
&P^{n} (x_i, \cdot)=\frac{\alpha}{2}\nu_1^i P^{n-n_1}+\frac{\alpha}{2}\(1-\frac{\alpha}{2}\)\nu_2^iP^{n-n_1-n_2}\\
&+\ldots+\frac{\alpha}{2}\(1-\frac{\alpha}{2}\)^{k-1}\nu_k^i P^{n-n_1-\ldots-n_k} +\(1-\frac{\alpha}{2}\)^{k}\mu_k^i P^{n-n_1-\ldots-n_k}
\end{align*}
for $i=1, 2$ and $n\ge n_1+\ldots+n_k$. From (3.9) we have
\begin{align*}
&\left|\int_X f(y)\nu_j^1 P^n (dy)-\int_X f(y)\nu_j^2 P^n (dy)\right|\\
&=\left|\int_X P^n f(y)\nu_j^1 (dy)-\int_X P^n f(y)\nu_j^2(dy)\right|\le\varepsilon/2
\end{align*}
for $j=1,\ldots, k$. By the definition of $k$ we then obtain
\begin{align*}
&|P^n f(x_1)-P^n f(x_2)|\\
&=\left|\int_X f(y) P^n (x_1, dy)-\int_X f(y) P^n (x_2, dy)\right|\\
&< \varepsilon /2 +2 \|f\|_{\infty} (1-\alpha /2)^k=\varepsilon.
\end{align*}
Since $\varepsilon>0$ and $f\in C_b(X)$ were arbitrary and since
linear combinations of point measures are dense in the space of all measures equipped with the weak topology, we have
$$
\operatornamewithlimits{lim}_{n\rightarrow\infty}\left|\int_X f(y)\mu_1P^n (dy)-\int_X f(y)\mu_2P^n (dy)\right|=0
$$
for all probability measures $\mu_1, \mu_2$ and for arbitrary $f\in C_b(X)$. Since $\mathbf \Phi$ is bounded in probability, the above condition is also satisfied for every $f\in C(X)$. On the other hand, from the above condition it follows that $\mu_*$ is a unique invariant measure and
$$
\mu P^n\stackrel{\rm w}{\rightarrow} \mu_*\qquad\text{as $n\rightarrow\infty$}
$$
for every probability measure $\mu$, which finishes the proof.
\end{proof}

As an immediate consequence of this theorem we obtain the following result due to \L. Stettner (see \cite{Stettner}):

\begin{cor} Assume that:\newline
{\bf (S1)} for every $\varepsilon >0$ and every compact set $K\subset X$ there exists a compact set $W\subset X$ such that\newline
\centerline{ $\inf_{x\in K} P^n (x, W)\ge 1-\varepsilon\qquad\text{for $n\in\mathbb N$,}$}\newline
{\bf (S2)} for every $f\in C_b(X)$ the functions $\{P^n f: n=1, 2,\ldots\}$ are equicontinuous on compact subsets of $X$,\newline
{\bf (S3)} for every open set $O\subset X$ and every $x\in X$\newline
\centerline{$P(x, O) >0,$}\newline
{\bf (S4)} there exist $\eta >0$ and a compact set $L\subset X$ such that for every compact set $W\subset X$\newline
\centerline{ $\inf_{x\in W} P^n (x, L)\ge \eta\qquad\text{for some $n\in\mathbb N$}.$}\newline
Then there exists a unique invariant measure $\mu_*$ for $\mathbf \Phi$ and $P^n (x, \cdot)$ convereges weakly to $\mu_*$.
\end{cor}

\section{A counterexample}

In the last section we shall define a discrete--time Markov--Feller chain 
which satisfies condition $(\mathcal E)$ but it has not an invariant measure.

Let $(\Omega, \mathcal F, \mathrm {Prob})$ be a probability space and let $\overline{\mathbb N}= \mathbb N\cup\{\infty\}$. Define $\mathbf x : \mathbb N\times \mathbb N\times \overline{\mathbb N}\to l^{\infty}$ by the following
$$
\mathbf x (i, j, k)=(i,\overbrace{\mathstrut 0, \ldots, 0}^{j-times}, 2^{-k}, \ldots).
$$
It is easy to see that $X=\mathbf x ( \mathbb N\times \mathbb N\times \overline{\mathbb N})$ is a closed subset of $l^{\infty}$. Consider the discrete--time Markov chain $\mathbf\Phi=(\Phi_n)_{n\ge 1}$ defined by the formula
$$
\Phi_n=\mathbf x (\zeta_n,\,\, \eta_n,\,\, \xi_n)\qquad\text{for $n\in\mathbb N$},
$$
where $\zeta_n, \eta_n\colon\Omega\to \mathbb N$ and $\xi_n\colon\Omega\to\overline{\mathbb N}$ are Markov chains satisfying
\begin{align*}
&\mathrm {Prob}\, (\zeta_{n+1}=i,\,\,\eta_{n+1}=j,\,\, \xi_{n+1}=k |\,\, \zeta_n=i_0,\,\,\eta_n=j_0,\,\, \xi_n=k_0)\\
&=
\left\{ \begin{array}{rll}
         p_1(i_0, k_0)& \mbox{for} & i=1,\,\,j=j_0+1,\,\, k=1; \\
         p_2(k_0)      & \mbox{for} & i=i_0,\,\, j=j_0+1,\,\, k=k_0+1; \\
         1-p_1(i_0, k_0)-p_2 (k_0) & \mbox{for} & i=i_0+1,\,\, j=j_0+1,\,\, k=k_0. \end{array}\right.
\end{align*}
Moreover, we assume that $p_2(k)= k^{-4}$ for $k\in\mathbb N$, $p_1 (i, k)=1-p_2 (k)$ for $k<i!$ and $p_1(i, k)=p_2 (k)$ for $k\ge i!$.
Further $p_1 (i, \infty)=p_2 (\infty)=0$. 

To show that $\mathbf\Phi$ satisfies Feller's property fix $f\in C(X)$ and $x_0\in X$. Let $x_n\to x_0$ as $n\to\infty$. Without loss of generality we may assume that $x_n=\mathbf x(i, j_n, k_n)$, $x_0=\mathbf x(i, 1,\infty)$ and $k_n\to\infty$ as $n\to\infty$. Then 
\begin{align*}
P f(x_n)&=p_1 (i, k_n) f(\mathbf x(1, j_n+1, 1))+p_2 (k_n) f(\mathbf x(i, j_n+1, k_n+1))\\
&+(1-p_1 (i, k_n) -p_2 (k_n)) f((i+1, j_n+1, k_n))\\
&\underset {n\to\infty}{\to} f((i+1, 1, \infty))=Pf(x_0).
\end{align*}

Now let $x=\mathbf x (i_0, j_0, k_0)$ be such that $k_0\neq \infty$. We will show that there exists $\vartheta>0$ such that
$$
P^n (x, U_0)\ge\vartheta\qquad\text{for $n\in\mathbb N$,}
\leqno(4.1)
$$
where $U_0=\{\mathbf x(i, j, k): i=k=1,\,\, j\in\mathbb N\}.$ Since $p_2 (k)=k^{-4}$ for $k\in\mathbb N$,
$p_1 (i, k)=1-p_2 (k)$ for $k<i!$ and $p_1(i, k)=p_2 (k)$ for $k\ge i!$,
we easily check that
$$
\sup_{n\in\mathbb N}\mathrm{E}[\xi_n|\,\,\zeta_0=i_0,\,\,\eta_0=j_0,\,\,\xi_0=k_0]<\infty.
$$
%have
%$$
%\mathrm  {Prob}\,(\sup_{n\in\mathbb N} \xi_n <\infty\,\,|\,\,\zeta_0=i_0,\,\,\eta_0=j_0,\,\,\xi_0=k_0)>0
%$$
%and hence 
Chebyshev's inequality shows now that there exists $M_0>i_0$ such that
$$
\inf_{n\in\mathbb N} \mathrm {Prob}\,( \xi_n \le M_0\,\,|\,\,\zeta_0=i_0,\,\,\eta_0=j_0,\,\,\xi_0=k_0)>0.
$$
From this and the fact that $p_1 (i, k)=1-p_2 (k)$ for $i<k!$ we obtain
$$
\gamma=\inf_{n\in\mathbb N} \mathrm {Prob}\, (\zeta_n\le M_0!,\,\,\xi_n \le M_0\,\,|\,\,\zeta_0=i_0,\,\,\eta_0=j_0,\,\,\xi_0=k_0)>0.
$$
By the Markov property we have
$$
P^{n} (x, U_0)\ge\gamma\cdot\min_{1\le i\le M_0!,\,\, 1\le k\le M_0}p_1 (i, k)\quad\text{for $n\in\mathbb N$,}
$$
which shows that condition (4.1) holds with
$$
\vartheta=\gamma\cdot\min_{1\le i\le M_0!,\,\, 1\le k\le M_0}p_1 (i, k).
$$

Let $z=(1,0,\ldots,)$. Fix an open set $U$ such that $z\in U$. Let $r>0$ be such that $B(z, r)\subset U$. Choose $k\in\mathbb N$ such that $x(1, j, k)\in B(z, r)$ for $j\in\mathbb N$. Then by the Markov property we obtain
$$
P^{n+k} (x, U)\ge\vartheta p_2(1)\cdot\ldots\cdot p_2 (k)\quad\text{for $n\in\mathbb N$,}
$$
which gives condition $(\mathcal E)$.

 Finally, it is obvious that $\mathbf\Phi$ does not admit an invariant measure since $\lim_{n\to\infty} \eta_n =\infty$ $\mathrm {Prob}$--a.s. \qed

\noindent {\footnotesize Institute of Mathematics,\\
Silesian University, Bankowa 14, 40-007 Katowice,\\
Poland, {\it e-mail: szarek@itl.pl}

\end{document}


\begin{thebibliography}{1}
\bibitem{Borovkov} 
\textsc{Borovkov, V.} (1991).
Ergodicity condition for Markov chains that are not connected with irreducibility in the sense of Harris. \textit{Sibirsk. Math. Zh.}
\textbf{32 (4)} 6--19.
\MR{1142065}
\end{thebibliography}
\end{document}

\end